\documentstyle[amssymb]{amsart}
\newcommand{\opname}[1]{\mathop{\fam0#1}}
\newcommand{\Vol}{\opname{Vol}\nolimits\ }
\newcommand{\st}{{\bigm|}}
\newcommand{\dual}{\sp\circ}
\newcommand{\R}{{\Bbb R}}
\newcommand{\E}[1]{\opname{E}\left[ #1 \right]}
\newcommand{\GL}{\opname{GL}}
\newcommand{\dx}{d \vec{x}}
\renewcommand{\d}{\partial}

\newtheorem{theorem}{Theorem}
\newtheorem{corollary}[theorem]{Corollary}
\begin{document}
\title{Another Low-Technology Estimate in Convex Geometry}
\author{Greg Kuperberg}
\address{Department of Mathematics \\ University of California \\
Davis, CA 95616}
\email{greg@@math.ucdavis.edu}
\thanks{The author was supported by an NSF Postdoctoral Fellowship,
        grant \#DMS-9107908.}
\subjclass{Primary 52A21; Secondary 46B03}
\date{December 4, 1994}
\begin{abstract}
We give a short argument that for some $C > 0$, every $n$-dimensional
Banach ball $K$ admits a $256$-round subquotient of dimension at
least $Cn/(\log n)$.  This is a weak version of Milman's quotient of
subspace theorem, which lacks the logarithmic factor.
\end{abstract}
\maketitle

Let $V$ be a finite-dimensional vector space over ${\Bbb R}$ and let $V^*$
denote the dual vector space.  A {\em symmetric convex body} or {\em (Banach)
ball} is a compact convex set with nonempty interior which is invariant under
under $x \mapsto -x$.  We define $K\dual \subset V^*$, the {\em dual} of a
ball $K \subset V$, by 
$$K\dual = \{y \in V^* \st y(K) \subset [-1,1]\}.$$
A ball $K$ is the unit ball of a unique Banach norm $||\cdot||_K$ defined by
$$||v||_K = \min \{t \st v \in tK\}.$$
A ball 
$K$ is an {\em ellipsoid} if $||\cdot||_K$ is an inner-product norm.
Note that all ellipsoids are equivalent under the action of $\GL(V)$.

If $V$ is not given with a volume form, then a volume such as
$\Vol K$ for $K \subset V$ is undefined.  However, some expressions
such as $(\Vol K)(\Vol K\dual)$ or $(\Vol K)/(\Vol K')$ for $K,K' \subset V$
are well-defined, because they are independent of the choice of a volume
form on $V$, or equivalently because they are invariant under $\GL(V)$
if a volume form is chosen.

An $n$-dimensional ball $K$ is {\em $r$-semiround} \cite{Szarek} if it
contains an ellipsoid $E$ such that
$$(\Vol K)/(\Vol E) \le r^n.$$
It is {\em $r$-round} if it contains an ellipsoid $E$ such that $K
\subseteq rE$. Santal\'o's inequality states that if $K$ is an
$n$-dimensional ball and $E$ is an $n$-dimensional ellipsoid,
$$(\Vol K)(\Vol K\dual) \le (\Vol E)(\Vol E\dual).$$
(Saint-Raymond \cite{Saint-Raymond}, Ball \cite{Ball}, and Meyer and
Pajor \cite{Meyer-Pajor} have given
elementary proofs of Santal\'o's inequality.)
It follows that if $K$ is $r$-round, then either $K$ or $K\dual$
is $\sqrt{r}$-semiround.

If $K$ is a ball in a vector space $V$ and $W$ is a subspace, we define
$W \cap K$ to be a {\em slice} of $K$ and the image of $K$ in $V/W$ to
be a {\em projection} of $K$; they are both balls.  Following Milman
\cite{Milman}, we define a {\em subquotient} of $K$ to be a slice of a
projection of $K$.  Note that a slice of a projection is also a
projection of a slice, so that we could also have called a subquotient
a proslice.  It follows that a subquotient of a subquotient is a
subquotient.  Note also that a slice of $K$ is dual to a projection of
$K\dual$, and therefore a subquotient of $K$ is dual to a proslice (or
a subquotient) of $K\dual$.

In this paper we prove the following theorem:

\begin{theorem} Suppose that $K$ is a $(2^{k+1}n)$-dimensional ball which is
$(2^{(3/2)^k}\cdot4)$-semiround, with $k \ge 0$.  Then $K$ has a
$256$-round, $n$-dimensional subquotient.
\label{thme}
\end{theorem}

\begin{corollary} There exists a constant $C > 0$ such that
every $n$-dimensional ball $K$ admits a $256$-round subquotient
of dimension at least $Cn/(\log n)$.
\end{corollary}

The corollary follows from the theorem of John that every $n$-dimensional
ball is $(\sqrt{n})$-round.

The corollary is a weak version of a celebrated result of Milman
\cite{Milman,Pisier}:

\begin{theorem}[Milman] For every $C > 1$, there exists $D > 0$, and for every
$D < 1$ there exists a $C$, such that every $n$-dimensional ball $K$ admits a
$C$-round subquotient of dimension at least $Dn$. \label{thmilman}
\end{theorem}

However, the argument given here for Theorem~\ref{thme} is simpler
than any known proof of Theorem~\ref{thmilman}.

Theorem~\ref{thmilman} has many consequences in the asymptotic
theory of convex bodies, among them a dual of Santalo's inequality:

\begin{theorem}[Bourgain,Milman] 
There exists a $C>0$ such that for every $n$ and for every $n$-dimensional
ball $K$,
$$(\Vol K)(\Vol K\dual) \ge C^n(\Vol E)(\Vol E\dual).$$ \label{thbm}
\end{theorem}

Theorem~\ref{thbm} is an asymptotic version of Mahler's conjecture,
which states that for fixed $n$, $(\Vol K)(\Vol K\dual)$ is minimized
for a cube. In a previous paper, the author \cite{Kuperberg:convex}
established a weak version of Theorem~\ref{thbm} also, namely that
$${\Vol(K)\Vol(K\dual) \ge (\log_2 n)^{-n}\Vol(E)\Vol(E\dual)} $$
for $n \ge 4$.  That result was the motivation for the present paper.

The author speculates that there are elementary arguments for both
Theorems~\ref{thmilman} and~\ref{thbm}, which moreoever would establish
reasonable values for the arbitrary constants in the statements of
these theorems.

\section{The proof}

The proof is a variation of a construction of Kashin \cite{Szarek}.  For
every $k$ let $\Omega_k$ be the volume of the unit ball in ${\Bbb R}^k$;
$\Omega_k$ is given by the formula
$${\pi^{k/2} \over \Gamma(\frac k2 + 1)}.$$
Let $V$ be an $n$-dimensional vector space with a distinguished
ellipsoid $E$, to be thought of as a round unit ball in $V$, so that
$V$ is isometric to standard ${\Bbb R}^n$ under $||\cdot||_E$.  Give
$V$ the standard volume structure $\dx$ on ${\Bbb R}^n$. In particular,
$\Vol E = \Omega_n$.  Endow $\d E$, the unit sphere, with the invariant
measure $\mu$ with total weight 1.  If $K$ is some other ball in $V$,
then
$$\Vol K  = \Omega_n\int_{\d E} ||x||_K^{-n} d\mu$$
and more generally,
$$\int_K ||x||_{E}^k \dx =
{n\Omega_n\over n+k}\int_{\d E} ||x||_K^{-n-k} d\mu.$$
Let $f$ be a continuous function on $\d E$.
Let $0<d<n$ be an integer and consider the space of $d$-dimensional
subspaces of $V$.
This space has a unique probability measure invariant under rotational
symmetry.  If $W$ is such a subspace chosen at random with respect
to this measure, then for any continuous function $f$,
\begin{equation}
\int_{\d E} f(x) d\mu = \E{\int_{\d (E \cap W)} f(x) d\mu}, \label{erandom}
\end{equation}
where $\mu$ denotes the invariant measure of total weight 1 on
$E \cap W$ also.
In particular, there must be some $W$ for which the integral of $f$ on the
right side of equation~\ref{erandom} is less than or equal to that of the
left side, which is the average value.

The theorem follows by induction from the case $k = 0$ and from the claim that
if $K$ is a $(2n)$-dimensional ball which is $r$-semiround, then $K$ has an
$n$-dimensional 
slice $K''$ such that either $K''$ or its dual is $(2r)^{2/3}$-semiround. In
both cases, we assume that $K$ is $r$-semiround and has dimension $2n$
and we proceed with a parallel analysis.

There exists an $(n+1)$-dimensional subspace
$V'$ of $V$ such that:
\begin{equation}
\int_{\d E'} ||x||_K^{-2n} d\mu \le {\Vol K \over
\Vol E} = r^{2n}, \label{evprime}
\end{equation}
where $E' = E \cap V'$.  Let $K' = V' \cap K$.  Then
\begin{equation}
\int_{\d E'} ||x||_K^{-2n} d\mu = {2n \over (n-1)
\Omega_{n+1}} \int_{K'}
||x||_{E'}^{n-1} \dx. \label{ekprime}
\end{equation}
Let $p$ be a point in $K'$ such that $s = ||p||_E$ is maximized; in particular
$K'$ is $s$-round Let $V''$ be the subspace of $V'$ perpendicular to $p$ and
define $K'' = V'' \cap K$ and $E'' = V'' \cap E$. The convex hull $S(K'')$ of
$K'' \cup \{p,-p\}$ is a double cone with base $K''$ (or suspension of
$K''$), and $S(K'') \subseteq K'$.  We establish an estimate that shows that
either $s$ or $\Vol K''$ is small.  Let $x_0$ be a coordinate
for $V'$ given by distance from $V''$.  Then
\begin{eqnarray}
\int_{K'} ||x||_{E'}^{n-1} \dx
& \ge & \int_{S(K'')} ||x||_{E'}^{n-1} \dx
> \int_{S(K'')} |x_0|^{n-1} \dx \nonumber
\\ & = & 2\int_0^s x_0^{n-1}(\Vol (1-{x_0 \over s})K'') dx_0 \nonumber
\\ & = & 2(\Vol K'')s^n \int_0^1 t^{n-1}(1-t)^n dt \nonumber
\\ & = & (\Vol K'')s^n {2(n-1)!n! \over (2n)!} \label{ekprimeprime}
\end{eqnarray}
We combine equations \ref{evprime}, \ref{ekprime}, and \ref{ekprimeprime}
with the inequality:
$$ {\Omega_n 4n(n-1)!n! \over \Omega_{n+1}(n-1)(2n)!} =
{2\Gamma({n + 3 \over 2})(n-2)!n! \over \sqrt{\pi} \Gamma({n+2 \over 2})(2n-1)!}
> 4^{-n} $$
(Proof:  Let $f(n)$ be the left side.  By Stirling's approximation,
$f(n)4^n \to 1$ as $n \to \infty$.  Since
$${f(n+2) \over f(n)} = {n^2 + 2n - 3 \over 4n^2 + 8n + 3} < \frac14,$$
the limit is approached from above.)  The final result is that
$${\Vol K'' \over \Vol E''} \le (2r)^{2n}s^{-n} $$

In the case $k = 0$, $r = 8$. Since $E'' \subseteq K''$,
$\Vol K'' \ge \Vol E''$, which implies that $s \le 4r^2 = 256$.
Since $K''$ is $s$-round, it is the desired subquotient of $K$.

If $k > 1$, then suppose first that $s \le (2r)^{4/3}$.  In this case
$K''$ is $(2r)^{4/3}$-round, which implies by Santal\'o's inequality
that either $K''$ or ${K''}\dual$ is $(2r)^{2/3}$-semiround.  On the
other hand, if $s \ge (2r)^{4/3}$, then $K''$ is
$(2r)^{2/3}$-semiround. In either case, the induction hypothesis is
satisfied.

\nocite{Bourgain-Milman}


\end{document}